\documentclass[reqno]{amsart}
\usepackage{amssymb}

\newtheorem{thm}{Theorem}
\newtheorem{lem}{Lemma}
\newtheorem{cor}{Corollary}
\theoremstyle{example}
\newtheorem{example}{Example}
\theoremstyle{remark}
\newtheorem{remark}{Remark}
\newcommand{\eop}{\hfill$\square$}
\renewcommand{\i}{\infty}

\newcommand{\ip}{\int_0^{\pi/2}}

\newcommand{\atn}{{\mathrm{tan}}^{-1}\,}
\newcommand{\asn}{{\mathrm{sin}}^{-1}\,}

\newcommand{\e}{\varepsilon}

\renewcommand{\th}{\theta}
\newcommand{\z}{\zeta}
\newcommand{\Z}{{\mathbf Z}}
\newcommand{\R}{{\mathbf R}}
\renewcommand{\(}{\left(}
\renewcommand{\)}{\right)}
\newcommand{\lb}{\displaystyle{\phantom\int}\!\!\!\!\!\!\!\left\{}
\newcommand{\rb}{\displaystyle{\phantom\int}\!\!\!\!\!\!\!\right\}}

\begin{document}
\title[Catalan's constant]{A class of series acceleration formulae for Catalan's
constant}
\author{David M.~Bradley}
\address{Department of Mathematics and Statistics\\
 University of Maine \\
        Orono, Maine 04469-5752\\ U.S.A.}
\email{bradley@math.umaine.edu; dbradley@member.ams.org}
\thanks{Research supported by the Natural Sciences and Engineering
Research Council of Canada (NSERC) and the Simon Fraser Center for
Experimental and Constructive Mathematics (CECM)}

\date{November 30, 1997}
\begin{abstract} In this note, we develop transformation formulae and
expansions for the log tangent integral, which are then used to
derive series acceleration formulae for certain values of Dirichlet
$L$-functions, such as Catalan's constant. The formulae are
characterized by the presence of an infinite series whose general
term consists of a linear recurrence damped by the central binomial
coefficient and a certain quadratic polynomial. Typically, the
series can be expressed in closed form as a rational linear
combination of Catalan's constant and $\pi$ times the logarithm of
an algebraic unit.
\end{abstract}

\maketitle

\section{Introduction}\label{sect1}
Catalan's constant may be defined by means of \cite{AS}
\begin{equation}\label{Gdef}
   G:= \sum_{k=0}^\i \frac{(-1)^k}{(2k+1)^2} = L(2,\chi_4),
\end{equation}
where $\chi_4$ is the non-principal Dirichlet character modulo 4.
It is currently unknown whether or not $G$ is rational.

The purpose of this note is to develop and classify acceleration
formulae for slowly convergent series such as~\eqref{Gdef}, based on
transformations of the log tangent integral. The simplest
acceleration formula of its type that we wish to consider is
\begin{equation}\label{G3}
   G
 =  \frac{\pi}{8}\log(2+\sqrt3)
   +\frac{3}{8}\sum_{k=0}^\i\frac1{(2k+1)^2\binom{2k}{k}},
\end{equation}
due to Ramanujan~\cite{BerndtI, Rama}. We shall see that Ramanujan's
formula~\eqref{G3} is the first of an infinite family of series
acceleration formulae for $G$, each of which is characterized by the
presence of an infinite series whose general term consists of a
linear recurrence damped by the summand in~\eqref{G3}.  In each
case, the series evaluates to a rational linear combination of $G$
and $\pi$ times the logarithm of an algebraic unit (i.e. an
invertible algebraic integer). Perhaps the most striking example of
this phenomenon is
\begin{equation}\label{G5}
   G
 = \frac{\pi}{8}\log\left(\frac{10+\sqrt{50-22\sqrt5}}
   {10-\sqrt{50-22\sqrt5}}\right)
   +\frac{5}{8}\sum_{k=0}^\i\frac{L(2k+1)}{(2k+1)^2\binom{2k}{k}},
\end{equation}
where $L(1)=1$, $L(2)=3$, and $L(n)=L(n-1)+L(n-2)$ for $n>2$
are the Lucas numbers \cite{Hardy}, (M2341 in \cite{Sloane}).

We shall see that series acceleration results such as~\eqref{G3}
and~\eqref{G5} have natural explanations when viewed as consequences
of transformation formulae for the log tangent integral, although we
should remark that Ramanujan apparently derived his
result~\eqref{G3} by quite different methods. The connection with
log tangent integrals is best explained by the equation
\begin{equation}\label{GFourier}
   G = -\int_0^{\pi/4}\log(\tan\th)\,d\th,
\end{equation}
obtained by expanding the integrand into its Fourier cosine series
and integrating term by term. It will be shown that Ramanujan's
result~\eqref{G3} arises from the transformation
\begin{equation}\label{T3}
   2\int_0^{\pi/4} \log(\tan\th)\,d\th
 = 3\int_0^{\pi/12}\log(\tan\th)\,d\th.
\end{equation}
The roccoco formula~\eqref{G5} arises in a similar manner from the
transformation
\begin{equation}\label{T5}
   2\int_0^{\pi/4}  \log(\tan\th)\,d\th
 = 5\int_0^{3\pi/20}\log(\tan\th)\,d\th
 - 5\int_0^{\pi/20} \log(\tan\th)\,d\th.
\end{equation}

Heuristically, one expects such transformations to succeed
because the reduced range of integration on the right hand
side, when re-expanded into a series, provides a continuous
analog of bunching together many terms of the original series.

\section{The Log Tangent Integral}\label{sect2}

There is a limitless supply of transformation formulae for the log
tangent integral.  In subsection 2.2, an infinite family of linear
relations, of which both~\eqref{T3} and~\eqref{T5} are members, will
be derived. These relations will be used in conjunction with the
series expansions given in subsection 2.1 to develop a corresponding
infinite family of series acceleration formulae which includes
both~\eqref{G3} and~\eqref{G5} as special cases.

\subsection{Series Expansions}
We shall be concerned with only two series expansions for the log
tangent integral.  These are given in Theorems~\ref{thm1}
and~\ref{thm2} below.
\begin{thm}\label{thm1}
For $0\le x\le \tfrac12\pi$,
\[
   \int_0^x\log(\tan\th)\,d\th
 = -\sum_{k=0}^\i \frac{\sin((4k+2)x)}{(2k+1)^2}.
\]
\end{thm}

\noindent{\bf Proof.} Expand the integrand into its Fourier cosine
series. Integrating term by term is justified by the fact that the
Fourier series is boundedly convergent on compact subintervals of
$(0,\tfrac12\pi].$ \eop

For us, the significance of Theorem~\ref{thm1} derives mostly from
the specialization $x=\tfrac14\pi$, which yields the
relationship~\eqref{GFourier} between Catalan's constant and the log
tangent integral. On the other hand, the expansion in powers of
sines provided by Theorem~\ref{thm2} below is more widely
applicable.

\begin{thm}\label{thm2}
For $0\le x\le \tfrac14\pi$,
\[
   \int_0^x \log(\tan\th)\,d\th
 = x\log(\tan x) - \frac14
   \sum_{k=0}^\i \frac{(2\sin 2x)^{2k+1}}{(2k+1)^2\binom{2k}{k}}.
\]
\end{thm}

\noindent{\bf Proof.} First integrate by parts, rescale, and use the
double angle formula for sine:
\begin{align*}
   \int_0^x \log(\tan\th)\,d\th - x\log(\tan x)
   & = -\int_0^x \frac{\th\sec^2\th}{\tan\th}\,d\th\\
   & = -\int_0^{2x}\frac{\th\,d\th}{4\tan(\tfrac12\th)\cos^2(\tfrac12\th)}\\
   & = -\int_0^{2x}\frac{\th\,d\th}{2\sin\th}\\
   & = -\int_0^{\sin(2x)}\frac{2t\,\asn t}{\sqrt{1-t^2}}\cdot\frac{dt}{4t^2}.
\end{align*}
Now employ the power series expansion~\cite{PiAGM}
\[
   \frac{2t\,\asn t}{\sqrt{1-t^2}}
 = \sum_{k=1}^\i \frac{(2t)^{2k}}{k\binom{2k}{k}},\qquad |t|<1,
\]
and integrate term by term. The result follows. \eop

In addition to Theorem~\ref{thm1}, the following representations
were also more or less known to Ramanujan, and can be easily
verified by differentiation:
\begin{align*}
   \int_0^x\log(\tan\th)\,d\th
&= x\log(\tan x) +
   \sum_{k=0}^\i\frac{(-1)^{k+1}(\tan x)^{2k+1}}{(2k+1)^2},
\quad 0\le x\le\tfrac14\pi,\\
   \int_0^x\log(\tan\th)\,d\th
&= (\tfrac12\pi - x)\log(\cos x) -
   \sum_{k=1}^\i \frac{(\cos x)^k(\sin kx)}{k^2},
\quad 0\le x\le\tfrac12\pi,\\
   \int_0^x\log(\tan\th)\,d\th
&= x\log(\tan x) +\tfrac12\pi\log(2\cos x)\\
& - \sum_{k=0}^\i\binom{2k}{k}
   \frac{(\cos x)^{2k+1}+(\sin x)^{2k+1}}{4^k (2k+1)^2},
\quad 0\le x\le\tfrac12\pi.\\
\end{align*}

\subsection{Transformation Formulae}
It will be convenient to define
\begin{equation}\label{Tdef}
   T(r)
   := \int_0^{r\pi}\log(\tan\th)\,d\th, \qquad 0\le r\le \tfrac12.
\end{equation}

Our development will provide two distinct transformation formulae
for the $T$-function:
the multiplication formula,
which expresses $T$ at odd multiples of the argument in terms of a
multitude of other $T$-values; and the reflection formula, which makes
it possible to restrict the domain to
the interval $0\le r\le\tfrac14$, and which will effect a number
of simplifications
in our intermediate calculations, as we shall see.

\begin{thm}\label{thm3}
For all $0\le r\le\tfrac12$,
the reflection formula
\[
   T(r) = T(\tfrac12-r)
\]
holds.
\end{thm}

\noindent{\bf Proof.} First, note that $T(\tfrac12)=0$.  This can be
seen either by putting $x=\tfrac12\pi$ in Theorem~\ref{thm1}, or by
observing that
\[
   T(\tfrac12)
 = \ip\log(\tan\th)\,d\th
 = \ip\log(\sin\th)\,d\th - \ip\log\sin(\tfrac12\pi-\th)\,d\th
 = 0.
\]
It follows that
\begin{align*}
   T(r)
 &= \int_0^{r\pi}\log(\tan\th)\,d\th
  = \ip\log(\tan\th)\,d\th - \int_{r\pi}^{\pi/2}\log(\tan\th)\,d\th\\
 &= \int_{\pi/2-r\pi}^0 \log(\tan(\tfrac12\pi - \th))\,d\th\\
 &= \int_{\pi/2-r\pi}^0 \log(\cot\th)\,d\th\\
 &= \int_0^{(1/2-r)\pi}\log(\tan\th)\,d\th\\
 &= T(\tfrac12-r),
\end{align*}
as stated. \eop

To prove the multiplication formula, we require the following
product expansion for the tangent function.

\begin{lem}\label{lem1}
Let $m=2n+1$ be an odd positive integer, and let $x\in \R$.
Then
\[
   \frac{\tan(mx)}{\tan(x)}
 = \prod_{j=1}^n \tan\(\frac{j\pi}{m}+x\)\tan\(\frac{j\pi}{m}-x\).
\]
\end{lem}

\noindent{\bf Proof.} Let $w=e^{ix}$.  Then
\begin{align*}
   \frac{\tan(mx)}{\tan(x)}
 &= \(\frac{w^{2m}-1}{w^{2m}+1}\)\(\frac{w^2+1}{w^2-1}\)\\
 &= \prod_{k=1}^n\(\frac{w^2-e^{2k\pi i/m}}{w^2-e^{(2k-1)\pi i/m}}\)
     \(\frac{w^2-e^{-2k\pi i/m}}{w^2-e^{-(2k-1)\pi i/m}}\)\\
 &= \prod_{k=1}^n\(\frac{we^{-k\pi i/m}-w^{-1}e^{k\pi i/m}}
     {we^{-(2k-1)\pi i/2m}-w^{-1}e^{(2k-1)\pi i/2m}}\)\\
 &\qquad\times\(\frac{we^{k\pi i/m}-w^{-1}e^{-k\pi i/m}}
     {we^{(2k-1)\pi i/2m}-w^{-1}e^{-(2k-1)\pi i/2m}}\)\\
 &= \prod_{k=1}^n\frac{\sin(k\pi/m-x)\sin(k\pi/m+x)}
     {\sin((2k-1)\pi/2m-x)\sin((2k-1)\pi/2m+x)}.\\
\end{align*}
After expressing the sines in the denominator in terms of cosines
and letting $j=n-k+1$, we have
\begin{align*}
   \frac{\tan(mx)}{\tan(x)}
 &= \prod_{j=1}^n \frac{\sin(j\pi/m-x)\sin(j\pi/m+x)}
     {\cos(j\pi/m-x)\cos(j\pi/m+x)}\\
 &= \prod_{j=1}^n \tan\(\frac{j\pi}{m}+x\)\tan\(\frac{j\pi}{m}-x\)
\end{align*}
as required. \eop

\begin{thm}\label{thm4}
Let $m=2n+1$ be an odd positive integer, and let $0\le r\le 1/(2m)$.
Then the multiplication formula
\[
   T(mr)
 = m \sum_{j=0}^n T\(\frac{j}{m}+r\)-m\sum_{j=1}^n T\(\frac{j}{m}-r\)
\]
holds.
\end{thm}

\noindent{\bf Proof.} By Lemma~\ref{lem1},
\begin{align*}
   T(mr)
 &= \int_0^{mr\pi}\log(\tan\th)\,d\th = m\int_0^{r\pi}\log(\tan(mx))\,dx\\
 &= m\int_0^{r\pi}\log(\tan x)\,dx
    +m\sum_{j=1}^n\int_0^{r\pi}\log\tan\(\frac{j\pi}{m}+x\)\,dx\\
 & +m\sum_{j=1}^n\int_0^{r\pi}\log\tan\(\frac{j\pi}{m}-x\)\,dx\\
 &= m T(r) + m\sum_{j=1}^n \lb T\(\frac{j}m+r\)-T\(\frac{j}m\)\rb\\
 & -m\sum_{j=1}^n \lb T\(\frac{j}m-r\) - T\(\frac{j}m\)\rb\\
 &= m\sum_{j=0}^n T\(\frac{j}{m}+r\)-m\sum_{j=1}^n T\(\frac{j}{m}-r\),
\end{align*}
as stated. \eop

To obtain transformations such as~\eqref{T3} and~\eqref{T5}, we
apply the reflection formula (Theorem~\ref{thm3}) and the
multiplication formula (Theorem~\ref{thm4}) with $r$ chosen so as to
express $T(\tfrac14)$ in terms of the $T$-function at values of the
argument less than $\tfrac14$.  The resulting transformations are
distinguished according to the parity of $n$ in the multiplier
$m=2n+1$.

\begin{thm}\label{thm5}
Let $n$ be an odd positive integer.
Then
\[
   G = -T\(\frac14\)
     = \frac{2n+1}{n+1}\sum_{j=1}^n (-1)^j T\(\frac{2j-1}{8n+4}\).
\]
\end{thm}

\noindent{\bf Proof.} Let $p$ be a nonnegative integer.  In the
multiplication formula, let $n=2p+1$, so that $m=4p+3$, and put
$r=1/(4m)$.  Then
\begin{align*}
   T\(\frac14\)
  &= m\sum_{j=0}^p \lb T\(\frac{4j+1}{4m}\)
   + T\(\frac{4(n-j)+1}{4m}\)\rb \\
& - m\sum_{j=1}^p\lb T\(\frac{4j-1}{4m}\)
   + T\(\frac{4(n-j+1)-1}{4m}\)\rb \\
& - m T\(\frac{4(p+1)-1}{4m}\).
\end{align*}

Applying the reflection formula (Theorem~\ref{thm3}) to each term in
the preceding sums yields the simplification
\[
   T\(\frac14\) = 2m\sum_{j=0}^p T\(\frac{4j+1}{4m}\)
        - 2m\sum_{j=1}^p T\(\frac{4j-1}{4m}\)
        - m T\(\frac14\).
\]
The preceding expression can be simplified further by combining
the two sums into a single alternating sum.  Thus,
\[
   -T\(\frac14\) = \frac{2m}{m+1}\sum_{j=1}^{2p+1} (-1)^j
           T\(\frac{2j-1}{4m}\).
\]
Writing $p$ and $m$ in terms of $n$ completes the proof. \eop

Theorem~\ref{thm6} below addresses the alternative case in which the
multiplier is congruent to 1 modulo 4.

\begin{thm}\label{thm6}
Let $n$ be an even positive integer.  Then
\[
   G
 = -T\(\frac14\)
 = \frac{2n+1}n\sum_{j=1}^n (-1)^{j+1} T\(\frac{2j-1}{8n+4}\).
\]
\end{thm}

\noindent{\bf Proof.} Let $p$ be a nonnegative integer.  In the
multiplication formula let $n=2p$, so that $m=4p+1$, and again put
$r=1/(4m)$. Then
\begin{align*}
   T\(\frac14\)
&= m\sum_{j=0}^{p-1}\lb T\(\frac{4j+1}{4m}\)
    + T\(\frac{4(n-j)+1}{4m}\)\rb\\
& - m\sum_{j=1}^p\lb T\(\frac{4j-1}{4m}\)+T\(\frac{4(n-j+1)-1}{4m}\)\rb\\
& + m T\(\frac{4p+1}{4m}\).
\end{align*}

Applying the reflection formula (Theorem~\ref{thm3}) to each term in
the preceding sums yields the simplification
\[
   T\(\frac14\)
   = 2m\sum_{j=0}^{p-1} T\(\frac{4j+1}{4m}\)
   - 2m\sum_{j=1}^p T\(\frac{4j-1}{4m}\)
   +m T\(\frac14\).
\]
The preceding expression can be simplified further by combining
the two sums into a single alternating sum.  Thus,
\[
   -T\(\frac14\)
   = \frac{2m}{m-1}\sum_{j=1}^{2p}(-1)^{j+1} T\(\frac{2j-1}{4m}\),
\]
Writing $p$ and $m$ in terms of $n$ completes the proof. \eop

\begin{example}\label{ex1}
Putting $n=1$ in Theorem~\ref{thm5} yields the transformation $
   2T(\tfrac14) = 3T(\tfrac1{12}),
$ which is a restatement of~(\ref{T3}). Putting $n=2$ in
Theorem~\ref{thm6} yields the transformation $
   2T(\tfrac14) = 5T(\tfrac3{20})-5T(\tfrac1{20}),
$
which is~(\ref{T5}).
\end{example}

\section{Applications to Series Acceleration}\label{sect3}
\subsection{Catalan's Constant}

\begin{thm}\label{thm7}
Let $n$ be an odd positive integer.  For nonnegative integers $k$,
define a sequence
\[
   F_n(k)
   := \sum_{j=1}^n \((-1)^{n-j+1}2\cos\(\frac{j\pi}{2n+1}\)\)^k,
\]
and let
\[
   u_n := \prod_{j=1}^n \(\tan\(\frac{2j-1}{8n+4}\)\pi\)^{(2j-1)(-1)^j}.
\]
Then $u_n$ is a unit algebraic integer, and Catalan's constant has
the series acceleration formula
\[
   G = \(\frac{\pi}{4n+4}\)\log u_n
      -\(\frac{2n+1}{4n+4}\)\sum_{k=0}^\i
       \frac{F_n(2k+1)}{(2k+1)^2\binom{2k}{k}}.
\]
\end{thm}

\noindent{\bf Proof.} Apply Theorems~\ref{thm1} and~\ref{thm2} to
the right hand side of Theorem~\ref{thm5}. Thus,
\begin{align}\label{Pf7a}
   G &= \(\frac{2n+1}{n+1}\)\sum_{j=1}^n (-1)^j
     \(\frac{2j-1}{8n+4}\)\pi\log\(\tan\(\frac{2j-1}{8n+4}\)\pi\)
     \nonumber\\
 & - \(\frac{2n+1}{4n+4}\)\sum_{j=1}^n (-1)^j
      \sum_{k=0}^\i \frac{\(2\sin((2j-1)\pi/(4n+2))\)^{2k+1}}
             {(2k+1)^2\binom{2k}{k}}\nonumber\\
 &= \frac{\pi}4\sum_{j=1}^n (-1)^j\(\frac{2j-1}{n+1}\)
    \log\(\tan\(\frac{2j-1}{8n+4}\)\pi\)\nonumber\\
 &  -\(\frac{2n+1}{4n+4}\)\sum_{k=0}^\i\frac1{(2k+1)^2\binom{2k}{k}}
    \sum_{j=1}^n (-1)^j \(2\sin\(\frac{2j-1}{4n+2}\)\pi\)^{2k+1}.
\end{align}

The inner sum in~\eqref{Pf7a} simplifies somewhat if the sines are
expressed in terms of cosines.  Thus,
\begin{align}
& \sum_{j=1}^n (-1)^j \(2\sin\(\frac{2j-1}{4n+2}\)\pi\)^{2k+1}\nonumber\\
&=\sum_{j=1}^n (-1)^j
\(2\cos\(\frac{2n+1-(2j-1)}{4n+2}\)\pi\)^{2k+1}\nonumber\\
&=\sum_{j=1}^n (-1)^j \(2\cos\(\frac{n-j+1}{2n+1}\)\pi\)^{2k+1}\nonumber\\
&=\sum_{j=1}^n (-1)^{n-j+1} \(2\cos\(\frac{j\pi}{2n+1}\)\)^{2k+1}.
   \label{Pf7b}
\end{align}
Substituting~\eqref{Pf7b} into~\eqref{Pf7a} completes the derivation
of the stated formula.

It now remains to show that $u_n$ is indeed
an algebraic unit (i.e. an invertible algebraic integer) as claimed.
Let $x=(2j-1)\pi/(8n+4)$.
Since the units in any ring form a multiplicative group, it suffices
to show that the numbers $t:=\tan x$ are all algebraic units, or
equivalently, that the numbers $t$ satisfy monic polynomials with integer
coefficients and constant term $\pm 1$.

From the addition formula for the tangent function, one sees that
$\tan(kx)$ is a rational function of $t$ for each nonnegative integer $k$.
Indeed, if polynomials $p_k, q_k\in \Z[t]$ are defined by the recursion
\begin{equation}\label{recur1}
   {p_{k+1}\choose q_{k+1}}
   = {\;1\;\;\; t\choose -t\;\; 1}{p_k\choose q_k}, \quad\mathrm{for}
   \quad k\ge 0,\quad {p_0\choose q_0}={0\choose 1},
\end{equation}
then for all nonnegative integers $k$,
\[
   \tan((k+1)x) = \frac{\tan x + \tan(kx)}{1-\tan(x)\tan(kx)}
   = \frac{tq_k+p_k}{q_k-tp_k} = \frac{p_{k+1}}{q_{k+1}}.
\]
Since $\tan((2n+1)x) = \tan((2j-1)\pi/4) = (-1)^{j+1}$, it follows
that $t = \tan((2j-1)\pi/(8n+4))$ satisfies the polynomial equation
\[
   p_{2n+1}(t)\pm q_{2n+1}(t) = 0.
\]
It remains to show that $p_{2n+1}\pm q_{2n+1}$ has both highest degree
coefficient and constant coefficient equal to $\pm 1$.

Let $k$ be an odd positive integer.  From the
recursion~\eqref{recur1}, it follows that
\begin{align}
   p_{k+2}+q_{k+2} &= (1-2t-t^2)p_k + (1+2t-t^2)q_k,\label{recur2}\\
   p_{k+2}-q_{k+2} &= (1+2t-t^2)p_k - (1-2t-t^2)q_k.\label{recur3}
\end{align}
An easy induction shows that the respective degrees of $p_k$ and
$q_k$ are $k$ and $k-1$, for all odd positive integers $k$.  This
fact, combined with a second induction, shows that the highest
degree coefficient of $p_k\pm q_k$ is equal to $\pm 1$ for all odd
positive integers $k$.  Finally,~\eqref{recur2} and~\eqref{recur3}
show that
\[
   p_{k+2}(0)\pm q_{k+2}(0) = p_k(0)\pm q_k(0)
\]
and so a final induction proves that $p_k\pm q_k$ has constant
coefficient equal to $\pm 1$ for all odd positive integers $k$. \eop

\begin{remark}\label{rem1}
Suppose $n$ is fixed, and we partition the algebraic numbers
\[
   (-1)^{n-j+1} 2\cos\(\frac{j\pi}{2n+1}\)
\]
into disjoint sets of mutual
conjugates.  Then the product of the minimum polynomials for
each set of conjugates is precisely the characteristic polynomial
of the linear recurrence satisfied by the sequence $\{F_n(k)\}_{k=0}^\i$.
\end{remark}

\begin{example}\label{ex2}
Putting $n=1$ in Theorem~\ref{thm7} gives
\[
   G = -\frac{\pi}8\log\(\tan\(\frac{\pi}{12}\)\)
       +\frac38\sum_{k=0}^\i\frac{(2\cos(\pi/3))^{2k+1}}{(2k+1)^2
       \binom{2k}{k}},
\]
which is Ramanujan's formula~\eqref{G3}.
\end{example}

Theorem~\ref{thm7} has its even counterpart in Theorem~\ref{thm8}
below.

\begin{thm}\label{thm8}
Let $n$ be an even positive integer.  For nonnegative integers $k$,
define a sequence
\[
   F_n(k) := \sum_{j=1}^n \((-1)^j 2\cos\(\frac{j\pi}{2n+1}\)\)^k,
\]
and let
\[
   u_n := \prod_{j=1}^n \(\tan\(\frac{2j-1}{8n+4}\)\pi\)^{(2j-1)(-1)^{j+1}}.
\]
Then $u_n$ is a unit algebraic integer, and Catalan's constant
has the series acceleration formula
\[
   G = \(\frac{\pi}{4n}\)\log u_n
   - \(\frac{2n+1}{4n}\)\sum_{k=0}^\i\frac{F_n(2k+1)}{(2k+1)^2\binom{2k}{k}}.
\]
\end{thm}

We omit the proof of Theorem~\ref{thm8}, as it closely mimicks the
proof of Theorem~\ref{thm7}. Instead, we derive the
formula~(\ref{G5}) which relates Catalan's constant and the Lucas
sequence.

\begin{cor}\label{cor1}
Let $L(1)=1$, $L(2)=3$, and $L(n)=L(n-1)+L(n-2)$ for $n>2$
be the Lucas numbers.  Then Catalan's constant has
the series acceleration formula
\[
   G
 = \frac{\pi}{8}\log\left(\frac{10+\sqrt{50-22\sqrt5}}
   {10-\sqrt{50-22\sqrt5}}\right)
 + \frac{5}{8}\sum_{k=0}^\i\frac{L(2k+1)}{(2k+1)^2\binom{2k}{k}}.
\]
\end{cor}

\noindent{\bf Proof.} Put $n=2$ in Theorem~\ref{thm8}.  Letting
$\phi := 2\cos(2\pi/5) = \tfrac12(\sqrt5-1)$ and $\tau :=
2\cos(\pi/5)  = \tfrac12(\sqrt5+1)$, we have
\[
   G
 =
   \frac{\pi}8\log\left(\frac{\tan(\pi/20)}{\tan^3(3\pi/20)}\right)
 - \frac{5}{8}\sum_{k=0}^\i\frac{\phi^{2k+1}-\tau^{2k+1}}
   {(2k+1)^2\binom{2k}{k}}.
\]
Now recall \cite{Hardy} that
\[
   L(k) = \(\frac{1+\sqrt5}2\)^k + \(\frac{1-\sqrt5}2\)^k
\]
for all nonnegative integers $k$.
It follows that
\[
   G
 = \frac{\pi}8\log\left(\frac{\tan(\pi/20)}{\tan^3(3\pi/20)}\right)
   +\frac58\sum_{k=0}^\i\frac{L(2k+1)}{(2k+1)^2\binom{2k}{k}},
\]
and so it remains only to verify the non-trivial denesting relationship
\begin{equation}\label{denest1}
   \frac{\tan(\pi/20)}{\tan^3(3\pi/20)}
 = \frac{10+\sqrt{50-22\sqrt5}}{10-\sqrt{50-22\sqrt5}}.
\end{equation}
To express the tangent values in~\eqref{denest1} in terms of
radicals, we follow \cite{Lewin}, p.~50.  Let $t:=\tan(\pi/20)$.
Then
\[
   \tan\frac{3\pi}{20} = \frac{3t-t^3}{1-3t^2}
   = \tan\(\frac{\pi}4 - \frac{2\pi}{20}\)
   = \frac{1-2t/(1-t^2)}{1+2t/(1-t^2)}.
\]
Equating the previous rational expressions in $t$ gives the quintic
equation
\[
   (t-1)^5 = 20 t^2(t-1),\quad \mathrm{or} \quad (t-1)^2 = 2t\sqrt5,
\]
since $t\ne1$.  Putting $t=(1-\e)/(1+\e)$, it follows that
$\e\sqrt{5+2\sqrt5}=\sqrt5$, and
\[
   \tan\frac{\pi}{20}
 = \frac{\sqrt{5+2\sqrt5}-\sqrt5}{\sqrt{5+2\sqrt5}+\sqrt5},   \quad
   \tan\frac{3\pi}{20} = \frac{\sqrt{5+2\sqrt5}-1}{\sqrt{5+2\sqrt5}+1}.
\]
Therefore, we may write
\begin{equation}\label{denest2}
   \frac{\tan(\pi/20)}{\tan^3(3\pi/20)}
 = \left(\frac{\sqrt{5+2\sqrt5}+1}{\sqrt{5+2\sqrt5}-1}\right)^3
   \frac{\sqrt{5+2\sqrt5}-\sqrt5}{\sqrt{5+2\sqrt5}+\sqrt5}
 = \frac{a+b}{a-b},
\end{equation}
where $a$ and $b$ are to be determined.
Cross multiplying and expanding both sides, we have
\begin{equation}\label{cross}
   5b(3+\sqrt5) = a(3-\sqrt5)\sqrt{5+2\sqrt5}.
\end{equation}
Since $(3-\sqrt5)/(3+\sqrt5) = \tfrac12(7-3\sqrt5)$,
we may write~(\ref{cross}) in the form
\[
   10 b = a \sqrt{(7-3\sqrt5)^2(5+2\sqrt5)} = a\sqrt{50-22\sqrt5}.
\]
Therefore, if in~\eqref{denest2}, we take $a=10$ and
$b=\sqrt{50-22\sqrt5}$, then~\eqref{denest1} holds, and the proof is
complete. \eop

\begin{remark}\label{rem2}
It is unlikely that~\eqref{denest1} will simplify any further.
Zippel \cite{Zippel} gives two formulae (caution: there are
misprints) for denesting expressions involving square roots. Borodin
et al~\cite{Borodin} show that these are the only two ways that such
expressions can be denested over the rational number field. In
particular, $\sqrt{50-22\sqrt5}$ cannot be denested, because
$50^2-5\times 22^2 = 80$ and $22^2\times 5^2 - 50^2\times 5 = -400$
are not squares of rational numbers.
\end{remark}

\subsection{Some Additional Examples}
One can derive additional acceleration formulae by specializing the
value of $x$ in Theorems~\ref{thm1} and~\ref{thm2} and equating the
two results. In general, convergence improves as the value of $x$
decreases.  The following selections provide a representative sample
of perhaps the most interesting results that can obtained using this
approach.

\begin{example}\label{ex3}
Putting $x=\tfrac14\pi$ gives (cf.~\eqref{Gdef})
\begin{equation}\label{slowest}
   G = L(2,\chi_4)
   = \frac12\sum_{k=0}^\i\frac{4^k}{(2k+1)^2\binom{2k}{k}},
\end{equation}
which Ramanujan~\cite{BerndtI} derived previously by other methods.
We remark that~\eqref{slowest} is actually a series
\emph{deceleration} result.  The reason for the poor convergence is
we have used the trivial transformation $T(\tfrac14)=T(\tfrac14)$
which fails to exploit the reduced range of integration present in
the other transformations.
\end{example}

\begin{example}\label{ex4}
Putting $x=\tfrac16\pi$ gives
\[
   L(2,\chi_6)
 = \frac{\pi\sqrt3}{18}\log 3
  + \frac12\sum_{k=0}^\infty\frac{3^k}{(2k+1)^2\binom{2k}{k}},
\]
where $\chi_6$ is the non-principal Dirichlet character modulo 6
(i.e. $\chi_6(5)=-1$).
\end{example}

\begin{example}\label{ex5}
Putting $x=\tfrac18\pi$ gives
\[
   L(2,\chi_8)
 = {{\pi\sqrt2}\over8}\log(1+\sqrt2)
 + \frac12\sum_{k=0}^\i\frac{2^k}{(2k+1)^2\binom{2k}{k}},
\]
where $\chi_8$ is the Dirichlet character modulo 8 given by $\chi_8(1)=
\chi_8(3)=1$, and $\chi_8(5)=\chi_8(7)=-1$.
\end{example}


\section*{Acknowledgements} I'm grateful to Chris Hill, Jonathan
Borwein, Petr Lison\v ek, and John Zucker for their helpful
observations.

\section*{Appendix} Here, we outline the role that inverse symbolic computation
-- in particular, Maple's integer relations algorithms -- played in
the discovery process.

A vector $\vec v=(v_1,v_2,\dots,v_n)$ of real numbers is said to
possess an integer relation if there exists a vector $\vec
a=(a_1,a_2,\dots,a_n)$ of integers not all zero such that the scalar
product vanishes, i.e. $a_1v_1+a_2v_2+\cdots+a_nv_n=0$.  In the past
two decades, several algorithms which recover $\vec a$ given $\vec
v$ have been discovered~\cite{Bailey1, Bailey2, Forcade, Hastad,
Lenstra}. One of these, ``LLL'' \cite{Lenstra}, has been implemented
in Maple V, and with its help, the authors of~\cite{Zeta}
and~\cite{Sigsam} discovered new formulae for values of the Riemann
Zeta function. The obstacle which initially confounded efforts to
extend the classical results
\[
   \z(2) = 3\sum_{k=1}^\i\frac1{k^2\binom{2k}{k}},\quad
   \z(3) = \frac52\sum_{k=1}^\i\frac{(-1)^{k-1}}{k^3\binom{2k}{k}},\quad
   \z(4) = \frac{36}{17}\sum_{k=1}^\i\frac1{k^4\binom{2k}{k}}
\]
to higher zeta values was circumvented
by the introduction of harmonic sums into the search space.
Thus, for example, by searching for an identity of the form
\[
   \z(7)
   = r_1 \sum_{k=1}^\i\frac{(-1)^{k+1}}{k^7 \binom{2k}{k}}
   + r_2 \sum_{k=1}^\i\frac{(-1)^{k+1}}{k^5 \binom{2k}{k}}
     \sum_{j=1}^{k-1}\frac{1}{j^2}
   + r_3 \sum_{k=1}^\i\frac{(-1)^{k+1}}{k^3 \binom{2k}{k}}
     \sum_{j=1}^{k-1}\frac{1}{j^4},
\]
we~\cite{Zeta} found
\[
   \z(7)
 = \frac{5}{2}\sum_{k=1}^\i\frac{(-1)^{k+1}}{k^7 \binom{2k}{k}}
 + \frac{25}{2}\sum_{k=1}^\i\frac{(-1)^{k+1}}{k^3 \binom{2k}{k}}
   \sum_{j=1}^{k-1}\frac{1}{j^4},
\]
and infinitely many more, as well as some lovely integral and
hypergeometric series evaluations, besides.

We suspected that a similar reverse-engineered approach might work
for certain Dirichlet $L$-series values, such as Catalan's constant,
but searching for similar variations on Ramanujan's
example~\eqref{G3} failed.  In view of the ornate complexity
of~\eqref{G5} and its relatives (Theorem~\ref{thm7} and
Theorem~\ref{thm8}), we can now understand the reason for this
failure. For a direct attack, one would have had to introduce, among
other things, logarithms of algebraic units into the model, so that
in effect, one would have needed to know beforehand the formula one
was searching for in order to find it. Models based on the inverse
tangent integral~\cite{Lewin, Rama} suffer the same drawbacks.   On
the other hand, the model based on the log tangent integral is
suited perfectly.

The author arrived at the log tangent integral model while
attempting to give an alternative proof of Ramanujan's acceleration
formula~\eqref{G3}.  It was found that the proof reduced to that of
proving the integral transformation~\eqref{T3}. Isolating the
$T$-function of section 3 for study was then a natural choice. After
directing Maple's integer relations finding algorithms to hunt for
linear relations amongst various $T$-values, the following list was
produced:
\begin{align}
   T(1/2) &= 0,\label{refl1}\\
   T(1/3) &= T(1/6),\label{refl2}\\
   T(1/8) &= T(3/8),\label{refl3}\\
   3T(4/9) &= T(1/3)+T(2/9)-3T(1/9),\label{mult9a}\\
   T(2/10) &= T(3/10),\label{refl4}\\
   T(1/10) &= T(2/5),\label{refl5}\\
   T(1/12) &= T(5/12),\label{refl6}\\
   2T(1/4) &= 3T(1/12),\label{mult3a}\\
   T(3/14) &= T(4/14),\label{refl7}\\
   T(5/14) &= T(1/7),\label{refl8}\\
   T(1/14) &= T(3/7),\label{refl9}\\
   3T(2/5) &= -3T(1/15)+T(1/5)+3T(4/15),\label{mult5a}\\
   3T(7/15) &= -3T(2/15)+3T(1/5)+T(2/5),\label{mult5b}\\
   15T(1/15) &= 15T(2/15)-5T(1/5)+9T(1/3)-10T(2/5),\label{tough}\\
   T(3/16) &= T(5/16),\label{refl10}\\
   T(1/16) &= T(7/16),\label{refl11}\\
   T(2/9) &= T(5/18),\label{refl12}\\
   T(1/9) &= T(7/18),\label{refl13}\\
   T(1/18) &= T(4/9),\label{refl14}\\
   3T(1/18)&= 3T(5/18)+T(1/3)-3T(7/18),\label{mult18}\\
   T(3/20) &= T(7/20),\label{refl15}\\
   T(1/20) &= T(9/20),\label{relf16}\\
   5T(3/20) &= 5T(1/20)+2T(1/4)=5T(7/20).\label{mult20}
\end{align}

Aside from trivial substitutions arising from the reflection formula
(Theorem~\ref{thm3}), the list evidently exhausts all linear
relations amongst $T$-values with rational arguments having
denominator no greater than $20$.  In fact, each list entry is a
consequence of the reflection formula and the multiplication formula
(Theorem~\ref{thm4}). For example,~\eqref{mult9a} follows from the
multiplication formula with $m=3$ and $r=1/9$.  The slightly
trickier~\eqref{tough} follows from three applications of the
multiplication formula.  One takes $m=3$ with $r=1/15$ and $r=2/15$,
and then one takes $m=5$ with $r=1/15$. This gives three equations.
Multiplying the first through by $5/2$, the second through by
$-5/2$, and the third through by $3/2$ and adding the three
resulting equations gives~\eqref{tough}.

From the list, it was easy to deduce and subsequently prove the
reflection formula.  At the same time, Chris Hill of the University
of Illinois used the $m=3$ case of Lemma~\ref{lem1} to
prove~\eqref{T3} i.e.~\eqref{mult3a}.  This broke the dam, leading
to the proof of Lemma~\ref{lem1}, the multiplication formula
(Theorem~\ref{thm4}), and the remaining results of
sections~\ref{sect2} and~\ref{sect3}.

\end{document}